\documentclass{statsoc}
\usepackage[a4paper]{geometry}
\usepackage{lscape,graphicx,bm}
\usepackage{epsfig,multirow,color}
\usepackage{ragged2e}
\usepackage{makecell,amsmath,amssymb,mathrsfs}
\usepackage{pdflscape}
\usepackage{rotating}
\usepackage{graphicx}
\usepackage{epstopdf}
\usepackage[]{threeparttable}
\usepackage{etoolbox,float}
\usepackage[caption = false]{subfig} 
\usepackage{enumerate}
\usepackage{comment}
\usepackage{natbib, float}

\makeatletter
\patchcmd{\@makecaption}
{\parbox}
{\advance\@tempdima-\fontdimen2} % decrease the width!
{}{}
\makeatother  

\newtheorem{theo}{Theorem}
\newtheorem{coro}{Corollary}[section]

\def\hat{\widehat}
\def\tilde{\widetilde}
\def\bfZ{\bm{Z}}

\def\bfz{\bm{z}}
\def\var{{\rm var}}
\def\bfV{\bm{V}_i}
\def\bfW{\bm{W}_i}
\def\scoreboot{T_{BS}}
\def\scorecal{T_{CS}}
\def\scorerobust{T_S}
\def\LRboot{T_{BL}}
\def\LRcal{T_{CL}}

\def\eqd{\stackrel{d}{=}}
\title[Robust Tests in Survival Analysis under Covariate-Adaptive Randomization]{Robust Tests for Treatment Effect in Survival Analysis 
	under Covariate-Adaptive Randomization}
\author{Ting Ye and Jun Shao}\address{Department of Statistics, University of Wisconsin-Madison}

\begin{document}

%\def\Box{\rule{2mm}{2mm}}
%\baselineskip = 8.5mm
%\parskip = 2.5mm

%\begin{quotation}
%\centerline
\begin{abstract}
Covariate-adaptive randomization is popular in clinical trials with sequentially arrived patients for balancing treatment assignments across prognostic factors which may have influence on the response. { However, existing theory on tests for treatment effect under covariate-adaptive randomization is limited to tests under linear or generalized linear models, although covariate-adaptive randomization has been used in survival analysis for a long time and its main application is in survival analysis.} Often times, practitioners would simply adopt a conventional test such as the log-rank test or score test to compare two treatments, which is controversial since tests derived under simple randomization may not be valid under other randomization schemes. In this article, we prove that the log-rank test valid under simple randomization is conservative in terms of type I error under covariate-adaptive randomization, and the robust score test developed under simple randomization is no longer robust under covariate-adaptive randomization. We then propose a calibration type log-rank or score test that is valid and robust under  both simple randomization and a large family of covariate-adaptive randomization schemes. Furthermore, we obtain Pitman's efficacy of log-rank and score tests to compare their asymptotic relative efficiency. Simulation studies about the type I error and power of various tests are presented under several popular randomization schemes. 
\end{abstract}


%\par

%\end{quotation}\par
%End---------------------------------------------------------------------------------------------------

%Introduction---------------------------------------------------------------------------------------------------

\section{Introduction}
In a clinical trial to compare the effectiveness of two treatments, simple randomization, also known as complete randomization, assigns patients independently into two treatment groups with equal probability. Simple randomization is widely accepted because it  provides a basis for statistical inference. 
In survival analysis, however, patients are not all available for simultaneous assignment of treatment but rather arrive sequentially and must be treated immediately. In such a case, simple randomization  may yield highly disparate  sample sizes between treatment arms across prognostic factors, e.g., institution, disease stage, prior treatment, sex, and age, which are known or thought to have significant influence on the response. Imbalance of treatment assignments across these factors may cause a confounding of the treatment effect and obscure the causation of active treatment itself to the observed effect. 

Many covariate-adaptive randomized treatment allocation schemes  have been proposed, which have advantages of minimizing imbalance between treatment groups across covariates or prognostic factors, reducing selection bias, minimizing accidental bias, and improving efficiency in inference \citep{Efron1971, Pocock1975, Wei1977, Wei1978b, Wei1978a, Weir2003}. One of the oldest covariate-adaptive randomization methods is the minimization procedure developed by \cite{Taves1974}. \cite{Pocock1975} generalized Taves' method and proposed to sequentially allocate patients with different probabilities to ensure treatment arms are marginally balanced within levels of each prognostic factor. A special case of Pocock and Simon's method is to apply the biased coin method \citep{Efron1971} to patients stratified by prognostic factors, which is referred to as the covariate-adaptive biased coin method.  Another popular method that has been extensively implemented in clinical studies is the permuted block design stratified  based on prognostic covariates \citep{Zelen1974}. A nice summary of the randomization schemes can be found in \cite{Schulz2018}.  As pointed out in \cite{Taves2010}, there are over 500 clinical trials which implemented Pocock and Simon's procedure to balance important covariates from 1989 to 2008. More recent examples of applying covariate-adaptive randomization can be found in \cite{Ploeg2010}, \cite{Fakhry2015}, \cite{Breugom2015}, \cite{Stott2017}, and \cite{Sun2018}.  The applications of covariate-adaptive randomization are not limited to clinical trials, and they are particularly relevant for randomized experiments with many interventions, for example, in mobile health.

In spite of the high prevalence of covariate-adaptive randomization, conventional tests are often utilized in practice (e.g., in 2018 {\em New England Journal of Medicine}, there are 9 articles using stratified permuted block design but  conventional  log-rank test or Wald's test), which has been raising concerns because statistical testing of the treatment effect should be performed using a test procedure valid under the particular randomization scheme used in data collection. Here, the validity of a test procedure refers to the type I error rate of the test is no larger than a given significance level, at least in the limiting sense.  \cite{Shao2010} initiated theoretical investigations on the validity of two sample $t$-test under covariate-adaptive biased coin randomization when the response follows a linear model. More results  based on $t$-test under linear or generalized linear models and other covariate-adaptive randomization methods are given by	\cite{Shao2013}, \cite{Ma:2015aa} and  \cite{Bugni:2017aa}. However, so far the study of testing hypotheses under covariate-adaptive randomization  in survival analysis is limited to empirical investigations, despite the fact that covariate-adaptive randomization has been used in survival analysis for a long time, and its main application is in survival analysis. The main reason for lack of theoretical results is that tests commonly used in survival analysis, such as the log-rank and score tests are highly non-linear, thus are more complicated to study than tests used under additive models. Moreover, censoring adds another layer of complexity. In this article, we will highlight the commonalities and differences between our results and the results derived under additive models.

To obtain a valid test, one approach is to construct a test under a correctly specified model that includes all covariates used in covariate-adaptive randomization \citep{Shao2010}. However, it is not always practical to include all covariates in the model for constructing a valid test procedure, mainly because it is difficult to incorporate some covariates in a model and the more covariates are included, the more likely the  model may be misspecified or too complicated to be useful.
For example, as described in \cite{Breugom2015}, 
a Cox proportional hazard model was adopted for survival data under 
a permuted block randomization of size six with stratification according to center, residual tumor, time between last irradiation and surgery, and preoperative treatment. 
It is difficult to specify a correct Cox model with these many covariates. 

Thus, it is of great interest to derive tests that are valid under covariate-adaptive randomization and robust against any model misspecification, either 
some or all covariates used in randomization are not included in the Cox
model, or the proportional hazard is completely wrong. 
The robustness against misspecification of Cox proportional hazard model 
has been studied in survival analysis 
\citep{Lin1989,KONG:1997aa,DiRienzo:2002aa}, but all results are limited to simple randomization.

 The purpose of our research is to establish a comprehensive theory for log-rank, score, or Wald's test of treatment effect in survival analysis under covariate-adaptive randomization. Our major contributions are three folds. First, we initiate the studies  beyond tests under additive models, where we develop novel technical tools to address the non-linear nature of tests in survival analysis. Second, the robust hypothesis testing results in survival analysis are generalized from simple randomization to a large family of covariate-adaptive randomization, where data are dependent. Third, we fill in the gap between theory and practice and provide some guidance in this field where covariate-adaptive randomization has its major application.
 
%  which provides  Our theory is based on a working Cox model and includes a study of the robustness of tests against misspecification of the Cox model, under covariate-adaptive randomization. {\blue Novel technical tools are developed to address the nonlinear nature of  survival data and the complexity of covariate-adaptive randomization.  }

After a detailed description of some popular covariate-adaptive randomization methods, in Section 2 we define the validity and conservativeness of tests, and the robustness. Section 3 presents asymptotic results of the score test under  covariate-adaptive randomization, which leads to results on the validity and conservativeness of score test. This theory also covers the popular log-rank test. Since the score test and log-rank test are conservative under covariate-adaptive randomization, we propose a calibration method to construct valid tests which are also robust against model misspecification. To compare powers of tests, results  on Pitman's asymptotic relative efficiency  are also included.  In Section 4, simulation studies are conducted to examine the finite sample performances of tests and our theory regarding type I error and power is supported by simulation results. Section 5 summarizes our main findings. All technical details are given in the Appendix in Supplemental Material. 

\section{Design and Hypothesis Testing} 
\label{sec:design}
\subsection{Design}
\label{subsec: design}
Let $n$ be the total number of patients in both arms, $\bfV$ be a vector of all measured and unmeasured covariates for the $i$th patient, which can be time-varying but the index $t$ is omitted,  
 $\bfZ_i$ be a subset of $\bfV$, a discrete time-independent baseline measured covariate with finitely many categories for which we want to balance by covariate-adaptive randomization, and  $I_i$ be the treatment indicator equaling $j$ if patient $i$ is assigned to treatment $j$, $j=0,1$. 

We first consider the design, i.e., the method of generating treatment assignments, $I_i$'s.
Treatment assignments under simple randomization are achieved by tossing a fair coin independently with 
the other variables so that $P(I_i=1)=P(I_i=0)=1/2$ for all $i$. To alleviate the imbalance between treatment groups, permuted block randomization is most frequently used, which can ensure the balancedness of two treatment groups at the end of every block. For example, with a block of $2b$ consecutively enrolled patients, 
exactly $b$ patients  are randomly allocated to one treatment. The block size $2b$ can remain fixed throughout the trial or varied under a pre-specified pattern. The permuted block randomization is very easy to implement and is quite effective in eliminating unbalanced design, but sometimes it is criticized to be too deterministic to result in selection bias. Biased coin \citep{Efron1971} is an alternative approach that can assure the imbalance to be controlled in probability without enforcing strict balancing. It assigns the $(k+1)$th patient according to \vspace{-1mm}
\begin{equation}
P(I_{k+1}=1)=\left\{\begin{array}{lc}
p,& D_{k}<0\\
1/2,& D_{k}=0\\
1-p, & D_{k}>0
\end{array},
\right.
\label{BCD}
\end{equation}
where $p$ is a prefixed probability greater than $1/2$, $D_0=0$ and $D_{k}=\sum_{i=1}^{k}(2I_i-1)$,  the difference between the number of patients in treatment 1 and treatment 0 after $k$ assignments have been made. The urn design \citep{Wei1977,Wei1978b,Wei1978a} belongs to the family of biased coin randomization with adaptively changing $p$. It assigns the $(k+1)$th patient according to (\ref{BCD}) but with $p$ being $p_k=\frac{1}{2}+\frac{\omega |D_{k}|}{2(2s +\omega k)}$, where $\omega, s$ are pre-specified non-negative real numbers. When $\omega=0$, it is the same as the simple randomization. Note that $p_k$ tends toward $\frac{1}{2}$ as $k$ increases with fixed $D_k$, indicating the urn design would force balancedness at the beginning of the treatment allocation, and approach simple randomization as the size of the trial increases.

The aforementioned three adaptive randomizations are not themselves covariate adaptive since they do not use any covariate information in treatment allocation. To balance across patients' prognostic profile, we could form strata by levels of $\bfZ_i$ and apply these randomization methods within each stratum. They are named stratified permuted block randomization, covariate-adaptive biased coin randomization, and stratified urn design, respectively. 

To characterize the property of covariate-adaptive randomization, we define the within strata imbalance as follows,

\begin{equation}
D_n(\bfz)= \sum_{i: \bfZ_i = \bfz }(2I_i-1) = \sum_{i: \bfZ_i = \bfz }I_i-\sum_{i: \bfZ_i = \bfz }(1-I_i).  \label{eq:D}
\end{equation}
The first simple property is %\vspace{-2mm}
\begin{description}
	\item	(D1)  $n_{\bfz}^{-1/2}D_n(\bfz)=o_p(1)$ for every $\bfz$, where $n_{\bfz}$ is the number of subjects with $\bfZ_i=\bfz$.%\vspace{-2mm}
\end{description}
Stratified permuted block randomization and covariate-adaptive biased coin randomization are examples of covariate-adaptive randomization satisfying (D1). Specifically, under stratified permuted block randomization, for every $\bfz$, $D_n(\bfz)$ is at maximum half of block size if the block size is fixed, or half of the last block size if the block size varies; thus $D_n(\bfz)$ is bounded. As for covariate-adaptive biased coin randomization, it is proved in \cite{Efron1971} that $D_n(\bfz)$ is bounded in probability for every $\bfz$. Another type of covariate adaptive randomization designs has the property %\vspace{-2mm}
\begin{description}
	\item
	(D2) $n_{\bfz}^{-1/2}D_n(\bfz)\xrightarrow{D}N(0,\nu_{{}_D})$ for every $\bfz$ and a $\nu_{{}_D}>0$, and $D_n(\bfz)$ and $D_n(\bfz')$ are independent for all $\bfz\neq\bfz'$.% \vspace{-2mm}	
\end{description}
Here, $\xrightarrow{D}$ denotes convergence in distribution as the sample size $n$ increases to infinity. Simple randomization and stratified urn design are examples of covariate-adaptive randomization satisfying (D2). Specifically, $\nu_{{}_D}=1$ under simple randomization and $\nu_{{}_D}=\frac{1}{3}$ under stratified urn design for any $s$ and $\omega$  \citep{Wei1978b,Wei1978a}.

Pocock and Simon's marginal method \citep{Pocock1975} is also widely used. It assigns patients using (\ref{BCD}) but with $D_{k}$ defined as a weighted sum of squared or absolute differences between number of patients over marginal levels of $\bfZ_i$. The covariate-adaptive biased coin is a special case of Pocock and Simon's method when $\bfZ_i$ is  one-dimensional. This method does not directly enforce balance in each stratum of $\bfZ_i $, but can be applied when the number of strata is too large and stratified randomization is infeasible. Unfortunately, neither (D1) nor (D2) holds under Pocock and Simon's marginal method, while \cite{Ma:2015aa} proved its marginal imbalance measure is bounded in probability. It is obvious that (D2) does not hold because treatment allocations are correlated across strata; an example illustrating (D1) does not hold is given in Section \ref{sec: simulation}. It is worth noting that \cite{Hu2012} modified Pocock and Simon's approach and proposed to use a balance measure that is a weighted sum of the overall imbalance, marginal imbalance and strata imbalance. Through carefully designing the weights, the  imbalance measure $D_n(\bfz )$ can be bounded in probability and (D1) holds.

Unless simple randomization is used,  $I_i$'s are dependent and each $I_i$ depends on the entire $\{\bfZ_1 ,..., \bfZ_n\}$. In this article, we focus on the balanced treatment allocation, i.e. $E(I_i|\bfZ_1,...,\bfZ_n )=1/2$, but our results can be easily extended to general cases. 

\subsection{ Hypothesis Testing}
In this subsection, we describe data collected under a given treatment assignment design and introduce some notation. 
	Let  $X^*_{ij}$ and 	$C_{ij}$ be the potential failure time and censoring time, respectively, for patient $i$ assigned to treatment $j$,  
$X_{ij} = \min (X^*_{ij}, C_{ij})$, $\delta_{ij}=1$ if $X_{ij} = X^*_{ij}$, and $\delta_{ij}=0$ if $X_{ij} = C_{ij}$. Let $\lambda^{(1)}(t,\bfV)$, $\lambda^{(0)}(t,\bfV)$ denote the true underlying hazard function of $X_{i1}^*$ and $X_{i0}^*$, respectively. For each patient, only one of the two treatments can be received, so the observed response with possible censoring for patient $i$ is $(X_i, \delta_i)$, where $X_i=I_iX_{i1}+(1-I_i)X_{i0}$ and   $\delta_i=I_i\delta_{i1}+(1-I_i)\delta_{i0}$. 
Throughout we assume that  $(X^*_{i0},C_{i0}, X^*_{i1}, C_{i1},\bfV)$, $i=1,... , n$, are independent and identically distributed (i.i.d.), and that the  following conditions hold.% \vspace{-3mm}
%	Let  $(X_{ij}, \delta_{ij})$ be the potential response if patient $i$ is assigned to treatment $j$, where $X_{ij}$ is the possibly right-censored failure time, $\delta_{ij}=1$ if $X_{ij}$ is the failure time and $\delta_{ij}=0$ if $X_{ij}$ is censoring time. For each patient, only one of the treatment can be received, so the observed response with possible censoring for patient $i$ is $(X_i, \delta_i)$, where $X_i=I_iX_{i1}+(1-I_i)X_{i0}$ and   $\delta_i=I_i\delta_{i1}+(1-I_i)\delta_{i0}$. Throughout we assume that  $(X_{i0},\delta_{i0}, X_{i1}, \delta_{i1},\bfV)$, $i=1,... , n$, are independent and identically distributed (i.i.d.), and that the  following conditions hold. \vspace{-3mm}
\begin{description}
	\item 
	(C1)  (randomization). $(X^*_{ij}, C_{ij}, \bfV)$'s and $I_i$'s  are  independent conditioned on $\bfZ_i$'s. \vspace{-3mm} 
	\item
	(C2) (non-informative censoring). $X^*_{ij}$ and $C_{ij}$ for patient $i$ are independent conditioned on covariate $\bfV$, for $j=0,1$.   \vspace{-3mm}   
\item  (C3) (treatment-independent censoring). $C_{i1}|\bfV\eqd C_{i0}|\bfV$, where $A \eqd B$ means that $A$ and $B$ are identically  distributed.
	% \vspace{-3mm}   
%	\item 
%(C4)  There exists a $\tau >0$, such that $P(X_{ij} \geq \tau )>0$ for both $j=0,1$.
		%There exists a neighbourhood ${\mathscr B}$ of $\beta_*$ such that for each $\tau<\infty$,
	%$$
		%\sup_{t\in [0,\tau],\beta\in\mathscr{B}} |S_n^{(0)}(\beta,t)-s^{(0)}(\beta,t)|=o_p(1),
		%$$
		%where $s^{(0)}(\beta,t)=E\{S_n^{(0)}(\beta,t)\}$, $o_p(1)$ denotes a quantity converges to 0 in probability.
%	\vspace{-2mm}
\end{description} 
Condition (C1) is  reasonable because (i) given $\bfZ_i$, $\bfV$ contains covariates not used in randomization, and (ii) treatment assignments do not affect the  potential failure time  and censoring, although they do affect the observed outcomes $X_i$'s and $\delta_i$'s through $X_i=I_iX_{i1}+(1-I_i)X_{i0}$ and   $\delta_i=I_i\delta_{i1}+(1-I_i)\delta_{i0}$. 	
All the  randomization designs described so far satisfy (C1).
Condition (C2) is typical  in survival studies. 
Condition (C3) is critical for the robustness property, because it guarantees that the score function defined later in (\ref{U}) has asymptotic mean zero under $H_0$, regardless of whether the model used to derive the score function is misspecified or not. A slightly weaker condition is also assumed by \cite{DiRienzo:2002aa} and \cite{KONG:1997aa}
in the case of simple randomization. This condition is reasonable under many realistic situations, but it is 
not fully general, since it requires that after adjusting for $\bfV$, the censoring distribution no longer depends on the treatment group.  When censoring is because of adverse effects, for example, it is related with the treatment group. But if the adverse events can be largely explained by patients' genotype or prognostic factors that are either measured or unmeasured, then (C3) can still be reasonable.

We are interested in testing whether there is a treatment effect, i.e., 
$$H_0: \lambda^{(1)}(t,\bfV)=\lambda^{(0)}(t,\bfV)\qquad \mbox{versus} \qquad H_1:\lambda^{(1)}(t,\bfV)\neq\lambda^{(0)}(t,\bfV)
$$
A test statistic ${T}$ is a function of observed data  constructed such that the null hypothesis $H_0$ is rejected if and only if $|T|>z_{\alpha/2}$, where $\alpha$ is a given significance level and $z_{\alpha/2}$ is the $(1-\alpha/2)$th quantile of the standard normal distribution. $T$ is said to be (asymptotically) valid if under $H_0$,
\begin{equation}
\lim_{n\rightarrow \infty} P(|T|>z_{\alpha/2})\leq \alpha
\label{valid}
\end{equation}
with equality holding for at least some parameter values under the null hypothesis $H_0$. $T$ is said to be (asymptotically) conservative if under $H_0$, there exists an $\alpha_0$ such that
\begin{equation}
\lim_{n\rightarrow \infty} P(|T|>z_{\alpha/2})\leq \alpha_0<\alpha. \label{conservative}
\end{equation}

Cox proportional hazard model is a very popular model in survival analysis. 
Suppose that tests of $H_0$ are based on fitting the following working Cox proportional hazard for the $i$th patient,
\begin{equation}
\lambda_0(t)\exp\{\theta I_i+\beta' \bfW\}, \label{Cox}
\end{equation} 
where $\bfW$ is an observed vector whose components are bounded functions of the components of $\bfV$, $\beta$ is an unknown parameter vector, $\beta'$ is the transpose of $\beta$, and $\lambda_0(t)$ is an unspecified baseline hazard function. The function in (\ref{Cox}) is a working hazard because it can be 
 unequal to the true hazard,  either 
 $\bfW$ may have arisen from mis-modeling and/or omitting components of $\bfV$, or the form of proportional hazard is not correct. Using (\ref{Cox}), we only need to observe $\bfW$, not necessarily the entire $\bfV$. 
 
 Under simple randomization,  the score and Wald tests based on hazard (\ref{Cox}) are asymptotically equivalent \citep{DiRienzo:2002aa}. We find that the same is true under covariate-adaptive randomization and, thus, we focus on the score test in the rest of this article. 
 
 With the working  hazard (\ref{Cox}), the partial likelihood function is 
 $$L(\theta,\beta)=\prod_{i=1}^n\left[\frac{\exp(\theta I_i+\beta'\bfW)}{\sum_{k=1}^{n}Y_k(X_i)\exp\left(\theta I_k+\beta'{\bm W}_k\right)}\right]^{\delta_i},$$
 where $Y_i(t)=I_iY_{i1}(t)+(1-I_i)Y_{i0}(t)$, 
 $Y_{ij}(t)=I(X_{ij}\geq t)$, and $I( \cdot )$ is the indicator function.
 The model-based score test for testing $H_0$ is 
 \begin{equation}
 T_M = n^{-1/2}U_\theta (0,\hat{\beta}_0)/\{\hat{A}(0,\hat{\beta}_0)\}^{1/2} , \label{scorem}\vspace{-1mm}
 \end{equation}
   where \vspace{-1mm}
 \begin{equation}\label{U}
 U_\theta(0,\beta)= \frac{\partial \log L(\theta , \beta )}{\partial \theta }\bigg|_{\theta = 0}  =\sum_{i=1}^{n}\int_0^{\tau}\left\{I_i-\frac{S_n^{(1)}(\beta,t)}{S_n^{(0)}(\beta,t)}\right\}dN_i(t),
 \end{equation}
 $\hat{A}(0,\beta)=-n^{-1}\partial^2\log L(\theta,\beta)/\partial\theta^2|_{\theta=0}
$, $S_n^{(r)}(\beta,t)=n^{-1}\sum_{i=1}^{n} Y_i(t)\exp\left\{\beta'\bfW\right\}I_i^r$,  $r=0,1$,   $N_i(t)=I_iN_{i1}(t)+(1-I_i)N_{i0}(t)$, $N_{ij}(t)=\delta_{ij}I(X_{ij}\leq t)$, 
 the upper limit $\tau$ in the integral is a point satisfying $P(X_{ij} \geq \tau )>0$ for $j=0,1$,
 and $\hat{\beta}_0$ is the maximum partial likelihood estimator of $\beta$ under constraint $\theta=0$. %i.e., $\hat{\beta}_0$ satisfies $U_\beta(0,\hat{\beta}_0)=0$. 
 From Theorem 2.1 of \cite{Struthers:1986aa}, $\hat{\beta}_0$ converges in probability to a unique  vector  $\beta_*$ under some regularity condition (see the conditions in Theorem 1 of our Section 3.1), regardless of whether  (\ref{Cox}) is a misspecified hazard or not. If the hazard in (\ref{Cox}) equals the true hazard function, then $\beta_*$ is the true value of $\beta$, and $T_M$ in (\ref{scorem}) is valid in the sense of (\ref{valid}) under simple randomization as well as under covariate-adaptive randomization, where the result for covariate-adaptive randomization
 follows directly from  the general result in \cite{Shao2010}.

However, the model-based score test $T_M$
is very fragile to model misspecification. 
 A test developed under working hazard  (\ref{Cox}) is said to be robust  if it is  valid according to (\ref{valid}) regardless of hazard misspecification. In what follows, robustness refers to  robustness  against misspecification of the true hazard function. 

Under simple randomization,  there are  discussions on constructing robust score and Wald tests  \citep{Lin1989,KONG:1997aa,DiRienzo:2002aa}. 
 In particular, \cite{Lin1989} proposed the following score test robust under simple randomization,  \vspace{-2mm}
\begin{equation}
\scorerobust=n^{-1/2}U_\theta(0,\hat{\beta}_0)/\{\hat{B}(0,\hat{\beta}_0)\}^{1/2},\vspace{-4mm}
\label{eq:robust score}
\end{equation}
where \vspace{-1mm}
\begin{equation}
\hat{B}(0,{\beta})=\frac{1}{4n}\sum_{i=1}^n \left[\delta_i-\sum_{j=1}^{n}\frac{\delta_jY_i(X_j)\exp\{{\beta}' \bfW\}}{nS_n^{(0)}({\beta},X_j)}\right]^2 . \label{varest}
\end{equation}
The difference between $T_M$ in (\ref{scorem}) and 
$T_S$ in (\ref{eq:robust score}) is the variance estimator in the denominator. 
The variance estimator in (\ref{varest}) is a robust estimator under simple randomization. 

 Besides tests based on hazard (\ref{Cox}), another popular test in survival analysis  is the log-rank test, which  does not use any covariate, 
 and is robust under simple randomization. 
% Note that when $\bfW$ in (\ref{Cox}) is zero for all  $i=1,...,n$, the model based score test is equivalent to the popular  which is also known to be asymptotically valid under simple randomization without 
	
 %The question we would like to address  is, 
%what is the behavior of log-rank test and score test  under covariate-adaptive randomization? Is the score test $T_S$ given in (\ref{eq:robust score}) robust under simple randomization still robust under covariate-adaptive randomization? 

\section{Theorem and Methods}
\vspace{-2mm}
\subsection{ Asymptotics of Score Test }
\label{subsec: robust score test}

	Under $H_0$, $\lambda^{(1)}(t,\bfV)=\lambda^{(0)}(t,\bfV)$. Thus, we  use $\lambda(t,\bfV)$ to denote the unspecified true hazard function of patient $i$ under $H_0$. Throughout this article, unless otherwise specified, the expectation $E$ is taken under $H_0$ with respective to the true  $\lambda(t,\bfV)$, not necessarily the hazard in  (\ref{Cox}).
	The numerator of score test can be expressed as
%\begin{eqnarray}
%n^{-1/2}U_\theta(0,\hat{\beta}_0)=n^{-1/2}U_\theta(0,\beta_*)+o_p(1),
%\end{eqnarray}
%where $o_p(1)$ denotes a quantity converges to 0 in probability.  It follows from \cite{KONG:1997aa} that under $H_0$, 
\begin{align}
n^{-1/2}U_\theta(0,\hat{\beta}_0)&=n^{-1/2}U_\theta(0,\beta_*)+o_p(1) \nonumber \\
&=n^{-1/2}\sum_{i=1}^{n}\int_0^\tau \left(I_i-\frac{1}{2}\right)\left\{dN_i(t)-p(t)Y_i(t)\exp(\beta_*'\bfW)dt\right\}+o_p(1)\nonumber\\
&=n^{-1/2}\sum_{i=1}^{n}\left\{I_iO_{i1}-(1-I_i)O_{i0}\right\}+o_p(1), \label{eq:U}
\end{align}
where the first two equalities are proved in the Appendix (Supplementary Material), $p(t)= {E}\left\{Y_i(t)\lambda(t, \bfV)\right\}/E\left\{Y_i(t)\exp\left(\beta_*'\bfW\right)\right\}$,  and\vspace{-1mm}
\begin{equation}
O_{ij}=2^{-1}\int_0^\tau \{dN_{ij}(t)- p(t)Y_{ij}(t)\exp(\beta_*'\bfW)dt\}. \label{oij}\vspace{-1mm}
\end{equation}

  Reformulating the score function as in (\ref{eq:U}) is a critical step. It helps to deal with the difficulty arising from the non-linearity of score function and the dependence. Under simple randomization, the sum in (\ref{eq:U}) has  i.i.d.\ terms, and thus the central limit theorem can be easily applied. 
 Under covariate-adaptive randomization, the   sum in (\ref{eq:U})  consists of dependent terms  due to the fact that $I_i$'s are dependent and each $I_i$ depends on  $ \{ \bfZ_i, i=1,...,n \}$, not just $\bfZ_i$.
To handle this problem under covariate-adaptive randomization, we consider the decomposition $n^{-1/2}U_\theta(0,\beta_*)=U_1+U_2+o_p(1)$, where 
\begin{align}	
U_1	&=n^{-1/2}\sum_{i=1}^{n}\left\{I_i(O_{i1}-E_i)-(1-I_i)(O_{i0}-E_i)\right\} ,\nonumber \\ 
U_2 & = n^{-1/2}\sum_{i=1}^{n}\left(2I_i-1\right)E_i =  n^{-1/2}\sum_{\bfz}D_n(\bfz)E\left(O_{ij}|\bfZ_i=\bfz\right) , \nonumber
\end{align}
$O_{ij}$ is given by (\ref{oij}), $D_n(\bfz)$ is given by (\ref{eq:D}), $E_i=E(O_{ij}|\bfZ_i)$, $j=0,1$.  Note that  (C2)-(C3) imply $E\left(O_{i1}|\bfV\right)=E\left(O_{i0}|\bfV\right)$ and consequently $E_i$ does not depend on $j$.
 
Because (C1) implies $E(U_1)=0$, and the condition $E(I_i|\bfZ_1,...,\bfZ_n )=1/2$ implies $E(U_2)=0$, $n^{-1/2}U_\theta(0,\beta_*)$ has asymptotic mean zero. Together with the central limit theorem applied to the conditional distribution of $U_1$ given $(I_i, \bfZ_i)$'s and the dominated convergence theorem, we show in the Appendix that
\begin{eqnarray}
U_1\xrightarrow{D} N\left(0,E\left\{\var (O_{ij}|\bfZ_i)\right\} \right),
\label{eq:u1 normality}
\end{eqnarray}
 Result (\ref{eq:u1 normality}) generally holds  for covariate-adaptive randomization methods discussed in Section 2, including the simple randomization. 
 
To derive the asymptotic distribution for $U_2$, we need some property of covariate-adaptive randomization, as discussed in Section \ref{subsec: design}. If (D1) holds, then it immediately  follows from its definition that 
$U_2=o_p(1)$ and hence 
\begin{eqnarray}\label{prop1}
n^{-1/2}U_\theta(0,\hat{\beta}_0) \xrightarrow{D} N\left(0,E\left\{\var (O_{ij}|\bfZ_i)\right\} \right).
\end{eqnarray}
 Under condition (D2), we have $U_2\xrightarrow{D} N\left(0,\nu_{{}_D}\var(E_i) \right) $. Also, we prove in the Appendix that $U_1$ and $U_2$ are uncorrelated, and hence
\begin{eqnarray}\label{prop2}
n^{-1/2}U_\theta(0,\hat{\beta}_0) \xrightarrow{D} N\left(0,E\left\{\var  (O_{ij}|\bfZ_i)\right\}+\nu_{{}_D}\var(E_i) \right)
\end{eqnarray}
Note that (\ref{prop1}) can be written as a special case of (\ref{prop2}) with $\nu_{{}_D} =0$. Formally, we have the following theorem.

\begin{theo}\label{theo:robust score test}
Let $\lambda(t,\bfV)$ be the true hazard  under $H_0$ and  (\ref{Cox}) be used as a working hazard. Assume (C1)-(C3) and that 
	$-n^{-1}\partial^2\log L(\theta,\beta)/\partial\beta \partial \beta' |_{\beta=\beta_*}$ converges in probability to a positive definite matrix. Under randomization with either (D1) or (D2),
	\begin{eqnarray}\label{U limit}
n^{-1/2}U_\theta(0,\hat{\beta}_0) \xrightarrow{D} N\left(0,E\left\{\var  (O_{ij}|\bfZ_i)\right\}+\nu_{{}_D}\var(E_i) \right),
	\end{eqnarray}
	where $\nu_{{}_D}=0$ when (D1) holds.
\end{theo}

Note that this theorem provides a unifying result that applies for both simple randomization and a large family of covariate-adaptive randomization. As explained in Section \ref{subsec: design}, simple randomization is characterized by property (D2) with $\nu_D=1$, while generally $\nu_D<1$ under covariate-adaptive randomization designs because they provide more balanced treatment assignments. Therefore, the score test $T_S$ (\ref{eq:robust score}) developed under simple randomization may not be robust under covariate-adaptive randomization with $\nu_D<1$. More specifically, we observe that the denominator of $T_S$,  i.e. $\hat{B}(0,\hat{\beta}_0)$ defined in (\ref{varest}) is  obtained by replacing unknown quantities in $n^{-1}\sum_{i=1}^n O_{i}^2$ with their empirical estimators, and $E(O_{ij})=0$, we conclude that under $H_0$, 
\begin{equation}
\hat{B}(0,\hat{\beta}_0) \xrightarrow{P} \var\left(O_{ij}\right) = E\left\{\var(O_{ij}|\bfZ_i)\right\}+ \var (E_i), \label{var}
\end{equation}
where $\xrightarrow{P}$ denotes convergence in probability. 
The following result shows the validity or conservativeness of $T_S$ in 
(\ref{eq:robust score}) under covariate-adaptive randomization.

\begin{coro}\label{theo:score test conservativeness}
	Under the same assumptions as in the Theorem \ref{theo:robust score test},  $T_S$ in (\ref{eq:robust score}) has the following property under $H_0$,  regardless of whether the hazard in (\ref{Cox}) is  misspecified or not:
	\begin{eqnarray}\label{limit1}
	\lim_{n\rightarrow\infty} P(|T_{S}|>z_{\alpha/2})=2\Phi\left(-z_{\alpha/2}\left[\frac{ E\left\{\var(O_{ij}|\bfZ_i)\right\}+ \var (E_i)}{E\left\{\var(O_{ij}|\bfZ_i)\right\}+\nu_{{}_D}\var\left(E_i\right)}\right]^{1/2}\right),
	\end{eqnarray}
	where $\Phi$ is the standard normal cumulative distribution function.
\end{coro}

A number of conclusions can be drawn from (\ref{limit1}). First, 
$T_S$ is valid if $\nu_{{}_D} =1$ (e.g., under simple randomization). Second, $T_S$ is also valid when $\var (E_i)=0$. Note that $\var(E_i)=0$ together with $E(E_i)=0$ means $E_i=0$ a.s., or equivalently 
\begin{equation}\label{eq:E}
\frac{E\left\{Y_{ij}(t)\lambda(t,\bfV)|\bfZ_i\right\}}{E\left\{Y_{ij}(t)\exp(\beta_*'\bfW)|\bfZ_i\right\}}=\frac{E\left\{Y_{ij}(t)\lambda(t,\bfV)\right\}}{E\left\{Y_{ij}(t)\exp(\beta_*'\bfW)\right\}}\  \ \mbox{a.s.},
\end{equation}
where $\lambda (t,\bfV)$ is the true hazard under $H_0$, not necessarily the working hazard in (\ref{Cox}) with $\theta =0$. 
A sufficient condition for (\ref{eq:E}) is that the  working hazard is the same as the true hazard.
Third, under covariate-adaptive randomization designs with $\nu_{{}_D} < 1$, $T_S$ is conservative when  $\var(E_i)>0$. Since $T_S$ is not valid unless (\ref{eq:E}) holds, which  almost requires correctness of the working model, the test $T_S$ robust under simple randomization is no longer robust under some popular covariate-adaptive randomization methods with $\nu_{{}_D} < 1$.

%\subsection{Conservativeness of log-rank test}
The aforementioned results on the score test can be applied to any kind of model misspecification. A special case is when $\bfW\equiv0$, $n^{-1/2}U_\theta(0,\beta)$ defined in (\ref{U}) equals the numerator of the popular log-rank test statistic. Denote $T_L$ as the log-rank test statistic,
\begin{equation}
{T}_{L}  =n^{-1/2} \sum_{i=1}^{n}\int_0^\tau\left\{I_i-\frac{\bar{Y}_1(t)}{\bar{Y}(t)}\right\}dN_i(t)
\Big/\hat{\sigma} , \label{Eq:logrank}
\end{equation}
where 
 $\bar{Y}_1(t)=\sum_{i=1}^nI_iY_i(t)$,  $\bar{Y}_0(t)=\sum_{i=1}^n(1-I_i)Y_i(t)$, $\bar{Y}(t)=\bar{Y}_0(t)+\bar{Y}_1(t)$, and 
 $$
 \hat{\sigma}^2  =n^{-1}\sum_{i=1}^n\int_0^\tau \frac{\bar{Y}_1(t)\bar{Y}_0(t)}{\bar{Y}^2(t)}dN_i(t). $$ 
 In addition, as proved in the Appendix, under $H_0$, the variance estimator 
\begin{equation}
\hat{\sigma}^2 \xrightarrow{P} \var\left(\tilde{O}_{ij}\right),%=\frac{1}{4}E(\delta_i)=\frac{1}{4}\int_0^\tau E\left\{Y_i(t)\lambda(t,0,\bfV)\right\}dt.  
\label{LR var}
\end{equation}  
We can derive the asymptotic distribution of log-rank test statistic $T_L$ using Theorem \ref{theo:robust score test}, simply by setting all the $\bfW$'s to zero. The following result implies the conservativeness of log-rank test $T_L$.
\begin{coro}\label{theo:log-rank conservativeness}
	Under the same assumptions as in the Theorem \ref{theo:robust score test},  $T_L$ in (\ref{Eq:logrank}) has the following property under $H_0$:
	\begin{eqnarray}\label{eq:log-rank conservativeness}
	\lim_{n\rightarrow\infty} P(|T_{L}|>z_{\alpha/2})=2\Phi\left(-z_{\alpha/2}\left[\frac{ E\left\{\var(\tilde{O}_{ij}|\bfZ_i)\right\}+ \var (\tilde{E}_i)}{E\left\{\var(\tilde{O}_{ij}|\bfZ_i)\right\}+\nu_{{}_D}\var\left(\tilde{E}_i\right)}\right]^{1/2}\right),
	\end{eqnarray}
	where $\Phi$ is the standard normal cumulative distribution function, $\tilde{E}_i = E(\tilde{O}_{ij}|\bfZ_i)$,\\$\tilde{O}_{ij}=2^{-1}\int_0^\tau \left\{dN_{ij}(t)-\tilde{p}(t)Y_{ij}(t)dt\right\}$ and $\tilde{p}(t)= E\left\{Y_i(t)\lambda(t, \bfV)\right\}/E\left\{Y_i(t)\right\}.$
\end{coro}

If $\nu_{{}_D} =1$ (e.g., simple randomization), the log-rank test is valid. 
If $\nu_{{}_D} <1$, the log-rank test is conservative unless
 $\tilde{E}_{i}=0$ a.s., which is the 
 unrealistic situation where  the covariate $\bfZ_i$ used for randomization is independent of the outcome. Therefore, we conclude that generally, the log-rank test is not robust but conservative under covariate-adaptive randomization.

\subsection{Constructing robust score tests}
\label{subsec:calibration}
The results in Section 3.1 tell us that the score test defined by (\ref{eq:robust score}), which is robust under simple randomization, is no longer robust to model misspecification under covariate-adaptive randomization with $\nu_{{}_D} <1$. The log-rank test, another robust test under simple randomization, is always conservative under covariate-adaptive randomization with $\nu_{{}_D} <1$. 

Results (\ref{prop2}) and (\ref{var})  reveal why the score test $T_S$ in (\ref{eq:robust score}) is not robust:
The variance estimator $\hat{B}(0,\hat\beta_0)$ does not take into account of the variability of $U_2$ reduced by covariate-adaptive randomization, and thus it may be  too large.
 Luckily $\var (E_i)=0$ when  the hazard in (\ref{Cox}) equals the true hazard functionso that $T_S$ is still valid,  but when the working hazard (\ref{Cox}) is  misspecified  and $\var (E_i) >0$, this variance estimator is too large and causes the conservativeness of $T_S$. The same can be said to the log-rank test $T_L$ in (\ref{Eq:logrank}), except that $T_L$ does not get the help from modelling so that $\hat\sigma$ in the denominator of (\ref{Eq:logrank}) is always too large under covariate-adaptive randomization with $\nu_{{}_D} <1$, which results in the conservativeness of log-rank test. 
 
  If we can construct a variance estimator for the numerator of (\ref{eq:robust score}) or 
(\ref{Eq:logrank}) that converges to the right variability under given randomization and $H_0$, then we can obtain a robust test, that is robust to model misspecification under a large family of randomization satisfying (D1) or (D2).

Under linear models, \cite{Shao2010} proposed a consistent variance estimator using the bootstrap method  which involves re-generating treatment indicators for every bootstrap dataset with the same randomization procedure applied in the original dataset. This method can preserve the randomization structure within every bootstrap sample, therefore, can intrinsically adapt to different randomization scheme. The price to pay is a large amount of computation, since treatment indicators have to be generated for every bootstrap sample. 

We now want to construct another variance estimator, which shares the advantages of bootstrap method and is computationally easy. From (\ref{prop2}) and Theorem 1, we just need to construct a consistent estimator of  
$E\left\{{\rm Var}(O_{ij}|\bfZ_i)\right\}+\nu_{{}_D}\var(E_i)$. Let
$\hat{O}_i=\frac{1}{2}\left\{\delta_i-\sum_{j=1}^{n}\frac{\delta_jY_i(X_j)\exp\{{\hat\beta}_0' \bfW\}}{nS_n^{(0)}({\hat\beta}_0,X_j)}\right\}$ and  $\hat{\sigma}^2_{\bfz}$ be the sample variance of $\hat{O}_i $'s within $\bfZ_i=\bfz$. Then, a consistent estimator of 
$E\left\{\var(O_{ij}|\bfZ_i)\right\}$ is
$n^{-1}\sum_{\bfz}n_{\bfz}\hat{\sigma}^2_{\bfz}$. Also, let $\hat{E}_{\bfz}$ be the sample mean of $\hat{O}_i$'s within $\bfZ_i=\bfz$. Because $E(E_i)=0$, a consistent estimator of $\nu_{{}_D}\var(E_i)$ is $n^{-1}\nu_{{}_D}\sum_{\bfz}n_{\bfz}\hat{E}_{\bfz}^2$. Note that $\nu_{{}_D}$ is a known constant for a given randomization. 
Therefore, we propose the following calibrated score test statistic
\begin{eqnarray}
\scorecal=U_\theta(0,\hat{\beta}_0)\Big/ \left\{\sum_{\bfz}n_{\bfz}\left(\hat{\sigma}^2_{\bfz}+\nu_{{}_D}\hat{E}^2_{\bfz}\right)\right\}^{1/2}.
\label{score cal}
\end{eqnarray}
The calibrated score test statistic $T_{CS}$ generalizes the robustness of $T_S$ to a large family of covariate-adaptive randomization, including simple randomization, under the conditions assumed in Theorem \ref{theo:robust score test}. It can also intrinsically adapt to randomization scheme of different balancing property by adjusting $\nu_D$. When $\nu_D=1$, i.e. under simple randomization, $T_{CS}$ degenerates to $T_S$.

For the log-rank test, since its numerator equals
$n^{-1/2}U_\theta(0,\beta)$ with $\bfW\equiv0$,
a calibrated log-rank test  $T_{CL}$ can be defined as 
(\ref{score cal}) with all $\bfW$'s replaced by $0$. 
This calibrated log-rank test is also robust under the conditions in Theorem 1. Again, with $\nu_D=1$, $T_{CL}$ degenerates to ordinary log-rank test $T_L$.

For applications, the proposed calibrated score test and calibrated log-rank test are robust to arbitrary model misspecification under covariate-adaptive biased coin randomization and stratified permuted block randomization, both of which satisfy (D1), and the stratified urn design and simple randomization, both satisfy (D2).  %Theoretical  properties of these tests under Pocock and Simon's marginal method are still unknown, due to the fact that Pocock and Simon's marginal method does not satisfy (D1) and whether it satisfies (D2) is unknown. 

\subsection{Asymptotic relative efficiency}

The reason we use  covariates and the working hazard (\ref{Cox}) is that it results in a more powerful test than the log-rank test without adjusting for covariates, when   (\ref{Cox}) correctly or nearly correctly specifies the true hazard function. The robustness is just an added guarantee that the test is still valid when working hazard (\ref{Cox}) is misspecified. %Under simple randomization, it is known that the score test or Wald test is asymptotically more powerful than the log-rank test, in terms of Pitman asymptotic relative efficiency. 

 In this subsection, we assume that working model $\lambda_0(t)\exp\{\theta I_i+\beta' \bfW\}$ is the true hazard function and study Pitman's asymptotic relative efficiency of the log-rank type and score type of test statistics. Note that when working model is the true model, the null hypothesis can be equivalently formulated as $H_0:\theta=0$, and we calculate Pitman's asymptotic relative efficiency under the contiguous alternative hypothesis $H_n:\theta_n=\gamma/\sqrt{n}$ with a constant $\gamma \neq 0$.    

Based on the general formula in \cite{KONG:1997aa} and our Theorem \ref{theo:robust score test}, under $H_n$ and any randomization we have discussed so far, we have 
$$ T_{CS} 
\xrightarrow{D}
N\left(\gamma\sigma_{S},1\right), 
$$
where $\sigma^2_S = E\left\{\var(O_{ij}|\bfZ_i)\right\}$. Pitman's efficacy of $T_S$ is then 
$\gamma^2 \sigma^2_S$. Therefore, covariate-adaptive randomization will not further boost efficiency once covariates are correctly adjusted through modelling.%This result also holds under $H_n$ and any covariate-adaptive randomization we discussed so far, because $E(O_{ij}|\bfW)=0$ a.s. which implies $E(O_{ij}|\bfZ_i)=0$ a.s. Furthermore, $T_S$ and $T_{CS}$ are asymptotically equivalent when (\ref{Cox}) is correct and hence $T_{CS}$ has the same result.

Under covariate-adaptive randomization, since the log-rank test $T_L$ is not valid but conservative, and it can be easily shown that $T_{CL}$ is uniformly more powerful than $T_L$, we now consider the robust test $T_{CL}$. It is shown in the Appendix that under $H_n$ and randomization with property (D1) or (D2), 
\begin{equation}
T_{CL} \xrightarrow{D}
N\left(\gamma\sigma_L^2/\sigma_{C},1\right),  \label{T_CL}
\end{equation} 
where 
\begin{eqnarray}
\sigma_{L}^2&=&\sigma^2_S -\frac{1}{4}\int_0^\tau \lambda_0(t)\Lambda_0(t)\left[E\left\{Y_i(t)e^{2\beta'\bfW}\right\}-\frac{\left[E\{Y_i(t)e^{\beta'\bfW}\}\right]^2}{E\left\{Y_i(t)\right\}}\right]dt\nonumber\\
\sigma^2_C &=&  E\{ \var (\tilde{O}_{ij}|\bfZ_i)\} + \nu_{{}_D} \var (\tilde{E}_i) 
= \sigma_S^2 - (1-\nu_{{}_D})\var(\tilde{E}_{i}).\nonumber
\end{eqnarray}
Thus, Pitman's efficacy of $T_{CL}$ is $\gamma^2\sigma^4_L/\sigma_C^2$. Comparing within $T_{CL}$, a straightforward observation is that the Pitman's efficacy of $T_{CL}$ decreases with increasing $\nu_D$. Since $0 \leq \nu_{{}_D} \leq 1$, $T_{CL}$ is more efficient under designs with $\nu_{{}_D} = 0$. Hence, the randomization itself can boost efficiency by achieving more balanced treatment allocation across `useful' covariates. The second observation is that under covariate-adaptive randomization based on $\bfZ$ and $\bfZ'$ given the same $\nu_D$, satisfying $\sigma(\bfZ)\subset \sigma(\bfZ')$ where $\sigma(\cdot)$ denotes the $\sigma$-field, $T_{CL}$ based on $\bfZ'$ will be more efficient. Therefore, utilizing more covariate information in the randomization procedure can increase efficiency.

%Since $T_{CL}$ is robust, it is better than $T_L$ under covariate-adaptive randomization. 
%The reason for this is because 
%covariate-adaptive randomization reduces the variance $\sigma_S^2$ to $\sigma_C^2$ by balancing treatment assignments and calibration leads to a correct variance estimator to avoid conservativeness, and it can be easily shown that $T_{CL}$ is uniformly more powerful than $T_L$ under covariate-adaptive randomization.

%Based on the aforementioned results, we conclude that when (\ref{Cox}) is correct,  $T_S$ and $T_{CS}$ (with any randomization) are more efficient than $T_L$ (with simple randomization).  Also,  $T_{CL}$ (with covariate-adaptive randomization) is  more efficient than $T_L$ (with simple randomization). 

 It remains to compare $T_{CS}$ and $T_{CL}$ under the same $\nu_D$. From previous discussions, we know that covariate-adaptive randomization can augment efficiency; on the other hand, adjusting for covariates through correct modelling can also increase efficiency. The following theorem compares these two approaches. 
 
 \begin{theo}
 	In addition to the assumptions stated in the Theorem \ref{theo:robust score test}, assume further that  the working hazard (\ref{Cox}) equals the true hazard function. Pitman asymptotic relative efficiency of $T_{CS}$ and $T_{CL}$ is 
\begin{eqnarray}
	ARE(T_{CS},T_{CL})=\sigma^4_L / (\sigma_{C}^2\sigma_S^2 ) \leq 1,
\end{eqnarray}
with equality holds if and only if $\beta=0$.
 \end{theo}
This theorem is proved in the Appendix. We therefore arrive at the conclusion that both correct modelling and covariate-adaptive randomization can boost efficiency, but the covariate information incorporated in the randomization cannot fully recover the efficiency loss due to not modelling. 

This is different from the result under  linear models as proved in \cite{Shao2010} that correct modelling and covariate-adaptive randomization that incorporates all the covariate information can achieve the same efficiency. The reason for this difference can be explained as follows. In linear models, the effects of modelling and covariate-adaptive randomization are the same. They both reduce the variance of the numerators of tests.  In survival analysis, due to the non-linear nature of the score tests (including log-rank test), correctly adjusting for covariates does not reduce the variance of the numerator. Instead, it increases the asymptotic limit of the numerator. Therefore, it is almost like these two approaches exert their effects of increasing efficiency through different pathways, thus they do not achieve the same effect.

\section{Simulation Results}
\label{sec: simulation}
In this section, two simulation studies are carried out to examine the Type I error and power of  tests under simple randomization and four most popular covariate-adaptive randomization methods, the covariate-adaptive biased coin randomization, stratified permuted block randomization, stratified urn design, and Pocock and Simon's marginal method. 

We know that in general, $D_n(\bfz)$ for Pocock and Simon's marginal method  does not have property (D1).  The following table provides a numerical evidence that ${\rm var}(D_n(\bfz ))/n$ does not tend to 0. 	
%	 but a theoretical proof is not available. 
%Whether Pocock and Simon's marginal method has property (D2) for some $\nu_{{}_D}$  is unfortunately unknown. 
In this section, we conduct simulation studies to empirically exam the performance of tests under Pocock and Simon's marginal method. 
\begin{table}
	\centering
	%	\caption{Variance of $D_n(\bfz )$ for Pocock and Simon's marginal method when 
	%		$\bfZ$ has independent binary components $Z^{(1)}$ and $Z^{(2)}$,  $P(Z^{(1)}=1)=0.5$, $P(Z^{(2)}=1)=0.3$, 		 simulated by 1000 Monte Carlo repetitions. }
	\label{table: strata imbalance}
	\begin{tabular}{cccccccccc}
		& \multicolumn{4}{c}{$\var(D_n(\bfz ))$} & & \multicolumn{4}{c}{$\var(D_n(\bfz ))/n$}     \\ 
		\cline{2-5} \cline{7-10} 
		
		$n$                  & $D_n(0,0)$ & $D_n(0,1)$ & $D_n(1,0)$ & $D_n(1,1)$ &                                  & $D_n(0,0)$ & $D_n(0,1)$ & $D_n(1,0)$ & $D_n(1,1)$               \\
		400                   & 22.76  & 25.25  & 22.00  & 25.24  &                       & 0.057                         & 0.063                 & 0.055                 & 0.063                 \\
		800                   & 46.88  & 47.56  & 45.39  & 47.15  &                       & 0.059                         & 0.059                 & 0.057                 & 0.059                 \\
		1200                  & 63.48  & 66.16  & 63.26  & 65.79  &                       & 0.053                         & 0.055                 & 0.053                 & 0.055                 \\
		1600                  & 91.82  & 91.75  & 92.75  & 94.56  &                       & 0.057                         & 0.057                 & 0.058                 & 0.059                 \\
		2000                  & 111.56 & 111.74 & 109.59 & 111.43 &                       & 0.056                         & 0.056                 & 0.055                 & 0.056               \\\hline 
	\end{tabular}
\end{table}

%Note that the covariate-adaptive biased coin randomization and stratified permuted block randomization satisfy the condition (a) in Theorem \ref{theo:score test conservativeness}; the stratified urn design satisfies the condition (b) in Theorem \ref{theo:score test conservativeness}; generally the Pocock and Simon's marginal method does not satisfy either one of the condition. 

The first simulation study considers the situation where (\ref{Cox}) correctly specifies the true hazard function. We consider the following four tests described in Section 3, the score test $T_S$ defined by (\ref{eq:robust score}), the log-rank test $T_L$ in (\ref{Eq:logrank}), 
the calibrated  score test $T_{CS}$ defined by (\ref{score cal}), and the 
calibrated log-rank test $T_{CL}$ described after (\ref{score cal}). 
We also include two bootstrap tests, $T_{BS}$ and $T_{BL}$, which are given by 
(\ref{eq:robust score}) and (\ref{score cal}), respectively, with the denominators replaced by the squared roots of bootstrap variance estimators as described in \cite{Shao2010}. 
The model-based  score test $T_M$ is omitted, because when (\ref{Cox}) correctly specifies the true hazard function, $T_M$ is asymptotically equivalent with $T_S$. Thus, we consider a total of six tests. 

The following three cases are considered, in which $U(a,b)$ denotes the uniform distribution on  interval $(a,b)$ and $C$ is the censoring time. %, $Exp(a)$ denotes the exponential distribution with rate $a$.
\begin{description}
	\item[Case 1.] The true hazard $=\lambda_0\exp\left(\theta I+1.5Z\right)$, $C\sim U(20,50)$, $Z$ is binary with $P(Z=1)=0.5$, and $Z$ is used in covariate-adaptive randomization. 
	\item[Case 2.] The true hazard $=\lambda_0\exp(\theta I+1.5Z_1-Z_{21}-0.5Z_{22})$, $C\sim U(20,40)$, $Z_1$ is binary with $P(Z_1=1) = 0.5$, $Z_2$ is discrete  with $P(Z_2=1)=0.4$, $P(Z_2=2)=0.3$, $P(Z_2=3)=0.3$,  $Z_{2k}$ is the indicator of $Z_2=k$, $Z_1$ and $Z_2$ are independent, and $Z_1$ and $Z_2$ are used in covariate-adaptive randomization.
	\item[Case 3.] The true hazard $=\lambda_0\exp(\theta I-1.5Z_1+0.5Z_2^2)$,
	 $C\sim U(10,40)$, $Z_1$ is binary with $P(Z_1=1)= 0.5$, $Z_2\sim N(0,1)$,   $Z_1$ and $Z_2$ are independent,	 and $Z_1$ and 
	  discretized $Z_2$ with $K$  equal probability categories  are used in 
	  covariate-adaptive randomization. 
\end{description}

Some quantities used in the simulation study are: $\lambda_0=\log 2/12$, 
%which means the median survival time is 12 months, 
the significance level  $\alpha=5\%$, the probability $p$ used in covariate-adaptive biased coin randomization and Pocock and Simon's marginal approach is $2/3$, the block size for stratified permuted block randomization is $2b=4$, the parameters used for stratified urn design is $s=1, \omega=1$, the sample size $n=200$ and $500$, and the bootstrap variance estimator is  approximated by Monte Carlo with size $200$.
%Note that $\LRcal$ and $\scorecal$ are well defined for randomization schemes satisfying condition (a) or (b), hence they are not suitable for Pocock and Simon's marginal method. For illustration purpose, $\LRcal$ and $\scorecal$ for Pocock and Simon's marginal method is the same as those for covariate-adaptive biased coin randomization and stratified permuted block randomization. 

The simulation Type I error based on 10,000 runs is shown in Table \ref{table:1st,type1} and the simulation power based on 2,000 runs is shown in Figures 1-3.
Because the bootstrap  and calibrated tests have similar performances, only the calibrated tests are presented in the figures. $\scorecal$ is also omitted from the figures because it is almost the same as $\scorerobust$ when the hazard in  (\ref{Cox}) equals the true hazard function.  For the ease of reading, 
in Figures 1-2, 
stratified permuted block design, stratified urn design, and simple randomization are to respectively represent $\nu_{{}_D}=0$, $\nu_{{}_D}<1$, and $\nu_{{}_D}=1$, while 
in Figure \ref{Fig:case3}, 
 only the stratified permuted block design is included in to represent  covariate-adaptive randomization methods.

 The following conclusions can be made from Table \ref{table:1st,type1} and
  Figures 1-3. \vspace{-2mm}
\begin{itemize}
	\item[(1)] The log-rank test $T_L$ is conservative under covariate-adaptive randomization. Because of this conservativeness, the power of $T_L$ under covariate-adaptive randomization is smaller than that under simple randomization when treatment effect is small, but the trend is reversed later when treatment effect is large.
	\item[(2)] The type I error of $\LRboot$, $\LRcal$, $\scoreboot$, and $\scorecal$ under covariate-adaptive biased coin randomization, stratified permuted block randomization, and stratified urn design are close to the nominal level 5\%, depicting the { robustness} of bootstrap and  calibrated  log-rank test and score test. %under condition (a) or (b) in Theorem \ref{theo:score test conservativeness}. 
	\item [(3)] Generally, the log rank type of tests $T_L$, $\LRboot$, and $\LRcal$ are not as powerful as the score type of tests $\scorerobust$, $\scoreboot$, and $\scorecal$ when model is correctly specified. This is different from the case of linear and generalized linear models  \citep{Shao2010,Shao2013}, where the bootstrap t-test and Wald test using the same covariate information have almost the same power. 
	\item [(4)] The log rank type of tests $\LRboot$ and $\LRcal$ are more powerful under designs satisfying (D1) than they are  under designs satisfying (D2) such as 
	the stratified urn design and simple randomization. 
	\item [(5)] The three score tests $\scorerobust,\scoreboot,$ and $\scorecal$ have almost the same power under different randomization schemes, when the model is correctly specified. 
	\item [(6)] When a discretized continuous covariate is used in covariate-adaptive randomization, $\LRboot$ and $\LRcal$ are more powerful when the  covariate is discretized into more categories. On the other hand, too many categories may cause sparsity of data in some  strata. 
\end{itemize}

Overall, the findings in simulation exactly coincides with asymptotic results given in Section 3.

The second simulation study aims to study the robustness and efficiency of tests when model is misspecified. The model-based score test $T_M$ is also included to evaluate its performance under model misspecification. Thus, there are a total of 7 tests. 
We consider the following three cases.
\begin{description}
	\item[Case 4.] The true hazard $=\lambda_0\exp(\theta I+Z_1-2Z_1Z_{21}+Z_1Z_{22})$, $C\sim U(20,50)$, where $Z_1$, $Z_2$, and $Z_{2k}$ are the same as those in Case 2. The  working hazard (\ref{Cox}) is
	$\lambda_0(t) \exp(\theta I+ \beta_1 Z_1+\beta_2 Z_2)$,  which is a misspecified hazard model since it does not include the interaction. $Z_1$ and $Z_2$ are used in covariate-adaptive randomization. 
	\item [Case 5.] The true hazard $=\lambda_0\exp(\theta I-0.5Z_1+1.5Z_2^2)$, $C | Z_1 \sim 10+{\cal E}(2Z_1)$, where $Z_1$ and $Z_2$ are the same as those in Case 3
	and ${\cal E}(a)$ denotes the exponential distribution with mean $a$. The covariate-adaptive randomization is carried out with $Z_1$ and discretized $Z_2$ with 4 levels. The working hazard (\ref{Cox}) is
	$\lambda_0(t) \exp(\theta I+ \beta_1 Z_1+\beta_2 Z_2)$
	without recognizing that the effect of $Z_2$ is quadratic.
	\item [Case 6.] The failure time $X^*$ does not follow Cox proportional hazard model, but $X^*=\exp(\theta I+1.5Z)+\epsilon$, where $Z\sim N(0,1)$ and $\epsilon\sim {\cal E}(1)$ are independent.
	 $C\sim U(10,20)$.   The covariate-adaptive randomization is carried out with $Z$ discretized with 4 levels.  The working hazard (\ref{Cox}) is $\lambda_0(t)\exp(\theta I +\beta Z)$, which is a misspecified hazard model.
\end{description}
The simulation results of  type I error 
are shown in Table \ref{table:mis,type1} and the simulation results of power for Case 6 are presented in Figure \ref{Fig:case6}.  Several conclusions can be obtained as follows.

\begin{itemize}
	\item[(1)] The type I error for $\LRboot$, $\LRcal$, $\scoreboot$ and $\scorecal$ are close to the nominal level 5\% under covariate adaptive biased coin randomization, stratified permuted block randomization and stratified urn design, indicating that the proposed tests are robust to any model misspecification.
	\item [(2)] $T_L$ and $\scorerobust$ may be conservative when model is misspecified. Although $\scorerobust$ is robust against model misspecification under simple randomization, it is not robust  under covariate adaptive randomization, which agrees with our asymptotic results. 
	\item [(3)] The calibration for log-rank  and  score tests does not apply to Pocock and Simon's marginal method because its asymptotic property is unfortunately unknown.
	 But empirical results show that bootstrapping can intrinsically adapt to different randomization schemes and have good performance even under Pocock and Simon's marginal method, although there is no theoretical confirmation for the bootstrap under the marginal method. 
	\item [(4)]  When model is misspecified, the model-based score test $T_M$ can be conservative, and can also have  inflated type I error. In other words, it is very fragile to model misspecification and can have unexpected performance.
	\item [(5)] The calibrated score test $T_{CS}$ based on a wrong  model can be less efficient than the calibrated log-rank test $T_{CL}$ without using any model.
\end{itemize}
\section{Conclusions}
\label{sec:discussion}
We derive a unified theory of robust hypothesis testing against model misspecification in simple randomization and a large family of covariate-adaptive randomization. Based on that, we further study asymptotic validity, conservativeness, and efficiency of log-rank and score tests for treatment effect with survival outcome when covariate-adaptive randomization is applied. Empirical results are included to complement our theory. Our results apply to simple randomization,
 covariate-adaptive biased coin randomization, stratified permuted block randomization. stratified urn design, and partially to Pocock and Simon's marginal method. The following are our main conclusions and recommendations. 
\begin{itemize}
	\item[(1)] 
	%Tests without correctly utilizing the covariates used in covariate-adaptive randomization are conservative, such as the log-rank test, or the score and Wald tests  robust under simple randomization. 
	The log-rank, score, and Wald test robust against model misspecification under simple randomization are not robust under covariate-adaptive randomization, but they are conservative. 
		% or Wald test omitting covariates used in stratification.
%	Regular robust score test developed under simple randomization can also be conservative under covariate-adaptive randomization when model is misspecified. 
	\item[(2)] The model-based score test or Wald test is valid only when working model equals the true model.  When model is misspecified, they may have inflated type I error. 
	\item[(3)] 
	A calibration is recommended for log-rank or score test that leads to robust tests.
	When the working model is true or nearly true, the calibrated  score test is more powerful than the calibrated log-rank test. However, their relative performance is unknown when the model is misspecified. 
\item[(4)]
 Covariate-adaptive randomization can boost efficiency by balancing treatment allocations across covariates. 
	Among different covariate-adaptive randomization methods, those can achieve better balancedness, such as 
 covariate-adaptive biased coin randomization and stratified permuted block randomization, 
	 lead to more powerful calibrated tests. Utilizing more covariate information in the design can also lead to more powerful calibrated tests.
	 \item[(5)]
	  Pocock and Simon's marginal method with score test works well if the working model equals the true model. When the working model is misspecified, however, its property is unknown, although its use together with the bootstrap performs well in empirical studies. 
	  \item[(6)] Under simple randomization, the log-rank test is popular because of its robustness against  model misspecification.  
	  The calibration technique developed in Section \ref{subsec:calibration} enhances its applicability to more general and better 
	  covariate-adaptive treatment randomization designs without sacrificing robustness.
\end{itemize}

The following is a summary of the performance of various tests under simple randomization (SR) or covariate-adaptive randomization (CA). Note that Wald's test is asymptotically equivalent to the corresponding score test. 
	\begin{table}
	\label{summary} 
		\centering
	\begin{tabular}{lcccccc}
		\hline
		& \multicolumn{5}{c}{When hazard is misspecified}&  \multicolumn{1}{c}{ Efficiency when (\ref{Cox})} \\
		\multicolumn{1}{c}{Method} & & under SR & & under CA & &   is the true hazard  \\ \hline
		model-base score test (\ref{scorem})& & incorrect${}^\dag$ & & incorrect${}^\dag$ && efficient \\
		score test (\ref{eq:robust score})& & robust & & conservative & & efficient \\
		log-rank test  (\ref{Eq:logrank})& & robust & & conservative & & inefficient \\
		calibrated score test (\ref{score cal})&& robust && robust && efficient \\
		calibrated log-rank test (\ref{score cal}), $\hat\beta_0=0$  && robust && robust && partially efficient \\
		\hline 
		\multicolumn{1}{l}{${}^\dag$ the type I error may be inflated} & \\
	\end{tabular}
	\end{table}

\section*{Supplemental Material}
The supplemental material contains the Appendix showing proofs of technical results in the article.

{
	\bibliographystyle{rss}
	\bibliography{reference}
}
 \begin{table}
	\caption{ 	\label{table:1st,type1} Simulation type I error in \% when  (\ref{Cox}) gives the true hazard. ($\alpha=5\%$, 10,000 runs)}
	\centering
		\resizebox{1\textwidth}{!}{
	\fbox{	\begin{tabular}{lccccccccccccc} \hline
			& \multicolumn{6}{c}{$n=200$}                           && \multicolumn{6}{c}{$n=500$}                           \\ \cline{2-7} \cline{9-14}
			& $T_L$ & $\LRboot$ & $\LRcal$ & $\scorerobust$ & $\scoreboot$ & $\scorecal$ && $T_L$ & $\LRboot$ & $\LRcal$ & $\scorerobust$ & $\scoreboot$ & $\scorecal$ \\
			{\footnotesize \bf Case 1} \\
			{\footnotesize Biased Coin}    &2.2   & 5.1 & 5.0   & 4.9   & 5.0   &   4.8&   & 1.7   & 4.8 & 4.6 & 4.7   & 4.9   & 4.6    \\
			{\footnotesize Permuted Block} & 2.0     & 4.8 & 4.7 & 4.5   & 4.6&   4.4&   & 2.2   & 5.4 & 5.1 & 5.2   & 5.2   & 5.1   \\
			{\footnotesize Marginal }      & 2.3   & 5.1 & - & 5.0     & 5.2   & 4.9&   & 1.9   & 5.1 &-  & 5.0     & 5.0     & 4.9   \\
			{\footnotesize Urn   }          &3.0     & 5.2 & 4.8 & 4.7   & 5.0     & 4.6   && 3.0     & 5.0   & 4.8 & 4.7   & 5.0     & 4.7   \\
			{\footnotesize SR   } &4.7   & 5.0   & 4.5 & 4.9   & 5.1   & 4.9  & & 4.6   & 4.8 & 4.5 & 4.8   & 4.9   & 4.8\\ %\cline{2-7} \cline{9-14}
			{\footnotesize \bf Case 2}            \\
			{\footnotesize Biased Coin}  & 1.9   & 5.3 & 5.8 & 5.0   & 5.2   & 4.9 &  & 1.6   & 5.0 & 5.0 & 4.9   & 5.0   & 4.8   \\
			{\footnotesize Permuted Block} &1.7   & 5.8 & 5.4 & 5.4   & 5.5   & 5.1   && 1.6   & 5.2 & 5.0 & 5.1   & 5.3   & 5.0   \\
			{\footnotesize Marginal }      &1.9   & 5.3 & - & 5.0   & 4.9   & 4.8   & &1.6   & 5.5 & -& 5.1   & 5.0   & 5.0   \\
			{\footnotesize Urn   }&2.7   & 5.3 & 5.0 & 5.6   & 5.7   & 5.4   & &2.6   & 5.2 & 5.0 & 5.3   & 5.4   & 5.2   \\
			{\footnotesize SR   }&4.7   & 4.9 & 4.5 & 4.3   & 4.6   & 4.3   && 5.0   & 5.1 & 5.0 & 5.0   & 5.3   & 5.0    \\ %\cline{2-7} \cline{9-14}
			{\footnotesize \bf Case 3, $K=8$}     &       &         &         &       &        &        &       &         &         &       &        &        \\
			{\footnotesize Biased Coin}    & 2.4   & 5.2 & 5.8 & 4.8   & 5.5   & 4.7 &  & 2.4   & 5.5 & 5.8 & 5.1   & 5.4   & 5.0   \\
			{\footnotesize Permuted Block} &2.0   & 6.2 & 5.5 & 5.4   & 6.1   & 5.3 &  & 1.7   & 5.5 & 5.2 & 5.3   & 5.8   & 5.2   \\
			{\footnotesize Marginal }      &2.4   & 5.5 & - & 5.0   & 5.4   & 4.9   & &2.0   & 5.1 & - & 4.8   & 5.0   & 4.8   \\
			{\footnotesize Urn   }&3.0   & 5.4 & 4.7 & 5.1   & 5.9   & 4.6   & &2.8   & 5.5 & 4.9 & 5.0   & 5.3   & 4.9   \\
			{\footnotesize SR   }&5.0   & 5.2 & 4.8 & 4.9   & 5.2   & 4.9   && 4.7   & 5.0 & 4.6 & 4.7   & 5.1   & 4.7   \\ %\cline{2-7} \cline{9-14}
			{\footnotesize \bf Case 3, $K=4$}	 &       &         &         &       &        &        &       &         &         &       &        &        \\
			{\footnotesize Biased Coin}    & 2.5   & 5.7 & 6.0 & 4.8   & 5.2   & 4.7 &  & 2.4   & 5.6 & 5.5 & 4.7   & 5.1   & 4.6   \\
			{\footnotesize Permuted Block} &2.2   & 5.7 & 5.3 & 5.1   & 5.6   & 5.0&   & 2.0   & 5.1 & 4.8 & 4.6   & 5.0   & 4.5   \\
			{\footnotesize Marginal }      &2.4   & 5.7 & - & 5.3   & 5.6   & 5.2   & &2.3   & 5.6 & - & 5.2   & 5.5   & 5.2   \\
			{\footnotesize Urn   }&3.7   & 6.0 & 5.6 & 5.1   & 5.4   & 4.9   & &3.0   & 5.3 & 4.9 & 4.6   & 4.7   & 4.5   \\
			{\footnotesize SR   }&5.0   & 5.2 & 4.8 & 4.9   & 5.2   & 4.9   & &4.7   & 5.0 & 4.6 & 4.7   & 5.1   & 4.7 \\ %\cline{2-7} \cline{9-14}
			{\footnotesize \bf Case 3, $K=2$}     &       &         &         &       &        &        &       &         &         &       &        &        \\
			{\footnotesize Biased Coin}  & 3.1   & 5.4 & 5.4 & 5.2   & 5.5   & 5.1 &  & 2.5   & 5.0 & 4.7 & 4.6   & 4.7   & 4.6   \\
			{\footnotesize Permuted Block} &2.6   & 5.0 & 4.7 & 4.6   & 4.8   & 4.5 &  & 2.6   & 5.3 & 4.9 & 5.4   & 5.5   & 5.4   \\
			{\footnotesize Marginal }&2.8   & 5.5 & -& 5.1   & 5.3   & 5.1   & &2.5   & 5.2 & - & 5.1   & 5.3   & 5.1   \\
			{\footnotesize Urn   }&3.6   & 5.4 & 5.1 & 5.0   & 5.3   & 4.9   & &3.4   & 5.3 & 5.0 & 5.2   & 5.5   & 5.2   \\
			{\footnotesize SR   }&5.0   & 5.2 & 4.8 & 4.9   & 5.2   & 4.9   & &4.7   & 5.0 & 4.6 & 4.7   & 5.1   & 4.7   \\ \hline
			\multicolumn{14}{l}{	$T_L$: log-rank test; $\LRboot$: bootstrap log-rank test; $\LRcal$: calibrated log-rank test}\\
			\multicolumn{14}{l}{ $\scorerobust$:  score test; $\scoreboot$: bootstrap score test; $\scorecal$: calibrated score test.}\\
			\multicolumn{14}{l}{Note: Under Pocock and Simon's marginal method, $T_{CL}$ is not applicable, while}\\
			\multicolumn{14}{l}{$T_{CS}$ is still valid because (\ref{Cox}) gives the true model.}\\
	\end{tabular}}}
\end{table}

\begin{table}
	\caption{\label{table:mis,type1}Empirical type I error in \% under two model misspecification cases. ($\alpha=5\%$, based on 10,000 runs).}
	\centering
	\resizebox{\textwidth}{!}{
		\fbox{	
			\begin{tabular}{lccccccccccccccc} \hline
				& \multicolumn{7}{c}{$n=200$}                           && \multicolumn{7}{c}{$n=500$}                           \\ \cline{2-8} \cline{10-16}
				&$T_M$& $T_L$ & $\LRboot$ & $\LRcal$ & $\scorerobust$ & $\scoreboot$ & $\scorecal$ &&$T_M$& $T_L$ & $\LRboot$ & $\LRcal$ & $\scorerobust$ & $\scoreboot$ & $\scorecal$ \\
				{ \footnotesize\bf Case 4} \\
				{\footnotesize Biased Coin}        & 3.2& 2.0   & 5.4 & 5.5 & 2.9   & 5.0   & 5.2&   &2.8& 1.8   & 4.8 & 4.8 & 2.6   & 4.8   & 4.7   \\
				{\footnotesize Permuted Block}    & 3.0  & 1.8   & 5.6 & 5.1 & 2.9   & 5.2   & 5.0 & & 3.0& 1.7   & 5.0 & 5.1 & 2.8   & 5.2   & 4.9   \\
				{\footnotesize Marginal }            & 5.6& 3.7   & 5.6 & - & 5.1   & 5.3   & -  & &5.8 &3.5   & 5.5 & - & 5.6   & 5.6   & -   \\
				{\footnotesize Urn   }                &4.0  & 2.7   & 5.4 & 4.9 & 3.7   & 5.3   & 5.1  && 3.6& 2.9   & 5.6 & 5.3 & 3.5   & 5.3   & 5.1   \\
				{\footnotesize SR   }                 &  5.3& 4.9   & 5.0 & 4.8 & 4.6   & 4.9   & 4.6  & &5.1& 5.1   & 5.0 & 5.0 & 4.8   & 5.0   & 4.8   \\
				&       &     &     &       &       &       &  &     &     &     &       &       &       \\
				{\footnotesize \bf Case 5} \\
				{\footnotesize Biased Coin}          & 2.3 & 2.3 & 5.5 & 5.8 & 2.3 & 5.2 & 5.7 &  &2.4& 2.3 & 5.8 & 5.8 & 2.4 & 5.8 & 5.8 \\
				{\footnotesize Permuted Block}     &  2.2& 2.2 & 6.0 & 5.5 & 2.2 & 5.7 & 5.4 &  &2.0& 1.9 & 5.2 & 5.0 & 2.0 & 5.3 & 5.0 \\
				{\footnotesize Marginal }            &2.5 & 2.3 & 5.3 & - & 2.4 & 5.3 & - &  &2.0& 1.9 & 5.0 & -& 2.0 & 5.0 & - \\
				{\footnotesize Urn   }                 & 2.7& 2.6 & 5.2 & 4.8 & 2.6 & 5.1 & 4.8 & & 3.1& 2.9 & 5.5 & 5.3 & 3.1 & 5.5 & 5.3 \\
				{\footnotesize SR   }                 &  4.7& 4.9 & 4.9 & 4.8 & 4.6 & 4.8 & 4.6 & & 4.9& 5.0 & 5.2 & 4.9 & 4.8 & 4.9 & 4.8  \\
				&       &     &     &       &       &       &  &     &     &     &       &       &       \\
				{\footnotesize \bf Case 6} \\
				{\footnotesize Biased Coin}&13.0 & 0.3 & 5.4 & 6.9 & 3.6 & 5.0 & 4.4& & 15.3 & 0.2 & 5.4 & 5.8 & 3.6 & 5.0 & 4.4 \\
				{\footnotesize Permuted Block}     &13.6 & 0.1 & 5.9 & 5.7 & 3.6 & 5.2 & 4.3 && 15.3 & 0.1 & 5.3 & 5.1 & 3.9 & 5.5 & 4.7\\
				{\footnotesize Marginal}     &13.7&	0.2	&5.0&	-&	3.5	&4.9	&-&&	14.8	&0.1	&5.1	&-	&3.9&	5.1&	-\\
				
				{\footnotesize Urn   }&14.0&	0.9&	5.4&	5.1&	3.4&	4.5	&3.8	&&15.7&	1.0&	5.6	&5.4	&4.1&	5.1&	4.6
				\\
				{\footnotesize SR   }  &14.6&	5.0	&5.2	&4.8&	4.1	&4.8&	4.1	&&15.6	&5.0&	5.2	&4.9&	4.2	&4.9&	4.2
				\\\hline
				\multicolumn{16}{l}{	$T_M$: model-based Wald test; $T_L$: log-rank test; $\LRboot$: bootstrap log-rank test; $\LRcal$: calibrated log-rank test}\\
				\multicolumn{16}{l}{ $\scorerobust$:  score test; $\scoreboot$: bootstrap score test; $\scorecal$: calibrated score test.}\\
				\multicolumn{16}{l}{Note: $T_{CL}$ and $T_{CS}$ are not applicable under Pocock and Simon's marginal method thus are omitted.}\\
	\end{tabular}}}
\end{table}

\begin{figure}
	\centering
	\makebox{\includegraphics[scale=0.55]{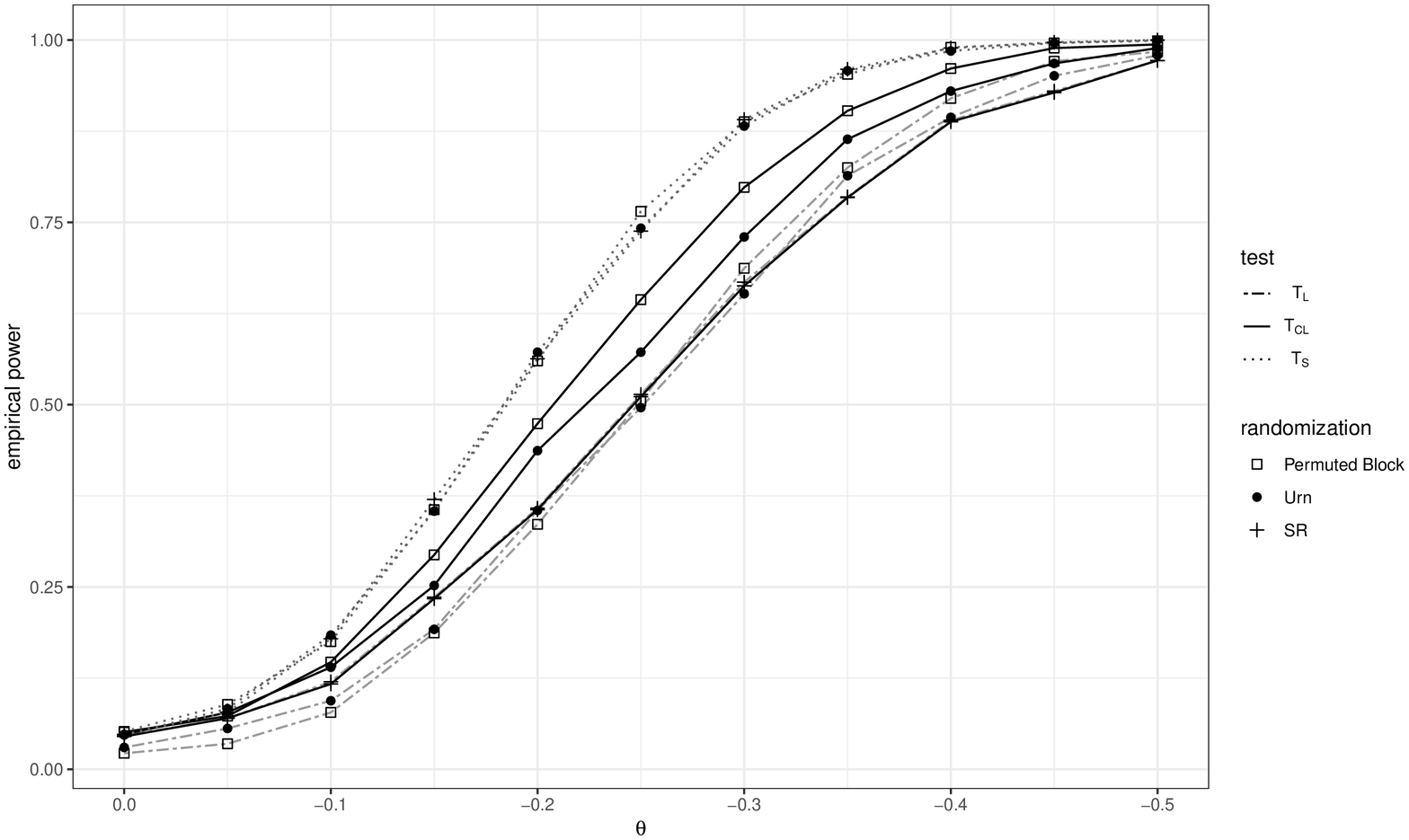}}
	\caption{\label{Fig:case1}Empirical power evaluated for $T_L$, $\LRcal$, $\scorerobust$ under three different randomization schemes in Case 1 ($\alpha=5\%$, $n=500$, based on 2000 runs).}
	\makebox{	\includegraphics[scale=0.55]{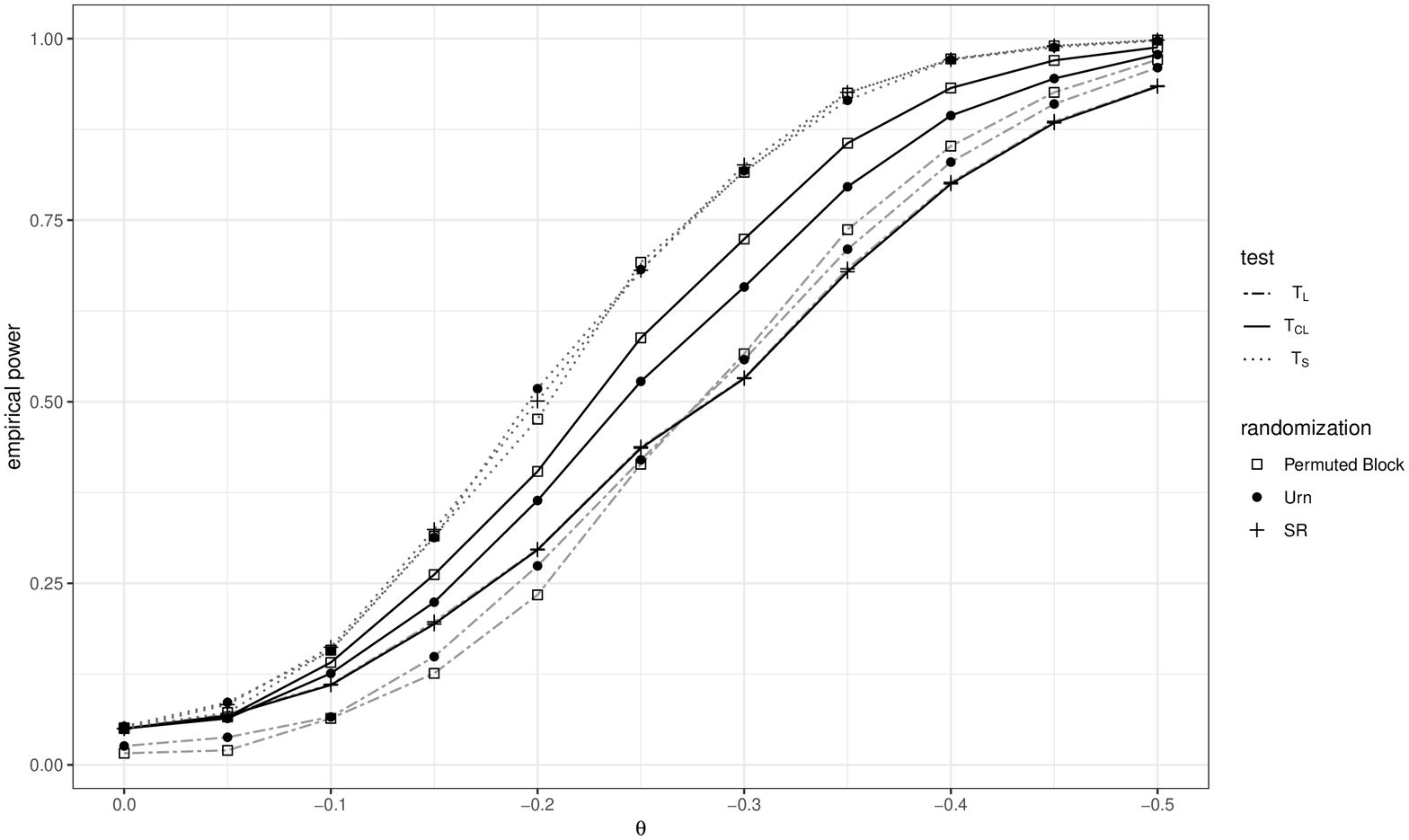}}
	\caption{\label{Fig:case1&2}Empirical power evaluated for $T_L$, $\LRcal$, $\scorerobust$ under three different randomization schemes in Case 2 ($\alpha=5\%$, $n=500$, based on 2000 runs)}
\end{figure}

\begin{figure}
	\centering
	\makebox{\includegraphics[scale=0.55]{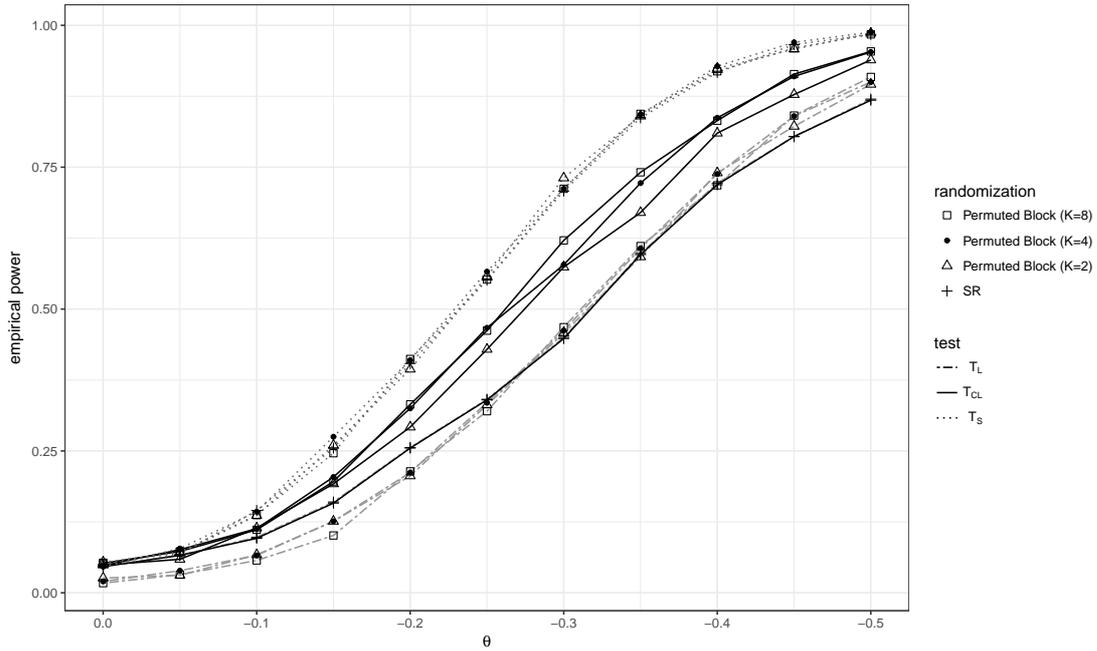}}
	\caption{	\label{Fig:case3}Empirical power evaluated for $T_L$, $\LRcal$, $\scorerobust$ in Case 3, with the continuous variable discretized into different number of categories ($\alpha=5\%$, $n=500$, based on 2000 runs)}
\end{figure}

\begin{figure}
	\centering
	\includegraphics[scale=0.55]{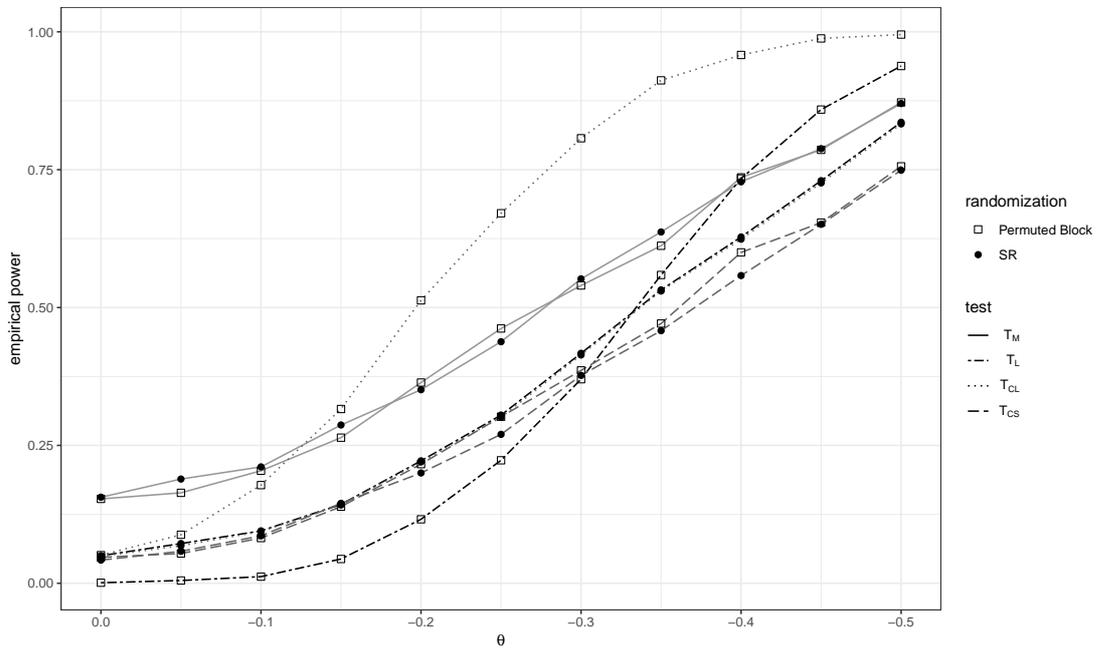}
	\caption{Empirical power evaluated for $T_M$, $T_L$, $\LRcal$, $\scorecal$ in Case 6. ($\alpha=5\%$, $n=500$, based on 2000 runs)}
	\label{Fig:case6}
\end{figure}

\end{document}